\documentclass[amssymb, 11pt]{amsart}
\usepackage{latexsym}
\newcommand{\C}{\mathcal{C}}
\newcommand{\F}{\mathbb {F}}

\newlength{\standardunitlength}
\newcommand{\aut}{{\mathrm {Aut}}}
\newtheorem{prop}{Proposition}[section]

\newtheorem{lemma}[prop]{Lemma}

\newtheorem{theorem}[prop]{Theorem}

\begin{document}

\title{Derangements in Simple and Primitive Groups}

\author{Jason Fulman}
\address{Department of Mathematics\\
University of Pittsburgh\\
Pittsburgh, PA 15260}
\email{fulman@math.pitt.edu}

\author{Robert Guralnick}
\address{University of Southern California \\ Los Angeles, CA}
\email{guralnic@math.usc.edu}



\thanks{Guralnick
was partially supported by National Science Foundation
grant DMS-9970305.}
\thanks{Originally submitted June 25, 2002, revised 
August 1, 2002}

\begin{abstract}  We investigate the proportion of fixed
point free permutations (derangements) in finite transitive
permutation groups.  This article is the first in a series where
we prove a conjecture of Shalev that the proportion of such
elements is bounded away from zero for a simple finite group.
In fact, there are much stronger results.
This article focuses on finite Chevalley groups of bounded
rank. We also discuss derangements in algebraic groups
and in more general primitive groups. 
These results have applications in questions about
probabilistic generation of finite simple groups and maps
between varieties over finite fields.
\end{abstract}

\maketitle

\section{Introduction}\label{intro}

Let $G$ be a group and $X$ a transitive $G$-set.  An element of $g \in
G$ is called a derangement on $X$ if $g$ has no fixed points on $X$.
We are interested in showing that under certain hypotheses the
set of derangements of $G$ on $X$ is large -- in particular, we will
mainly focus on the case where $G$ is finite.  We then define
$\delta(G,X)$ to be the proportion of elements in $G$ that are
derangements acting on $X$.  The rare situations when $\delta(G,X)$
is very small are also quite interesting and arise in
the theory of permutation and exceptional polynomials, coverings
of curves and  graph theory.  

The study of derangements goes back to the origins of permutation
group theory. It is an elementary result of Jordan that if $X$ is a
finite transitive $G$-set of cardinality $n > 1$, then
$\delta(G,X) > 0$.  It is also one of the earliest
problems in probability theory -- the problem was considered
by Monmort in 1708.  Diaconis pointed out to us that Frobenius in 1904
showed that $G \le S_n$ is $k$-fold transitive if and only if the first
$k$ moments of the number of fixed points is equal to the first $k$ moments of
a Poisson(1) random variable. He used this to determine character tables
of Mathieu groups.  

Jordan's result fails if $G$ is infinite.
There are various constructions for example
where any two nontrivial elements of $G$ are conjugate.  Then $G$
contains no derangements in any transitive action with a nontrivial
point stabilizer. Another
example is the case of $GL(V)$ where $V$ is a finite dimensional
vector space over an algebraically closed field  and $X$ is the
set of subspaces of fixed dimension or more generally $X$ is the
set of flags of a given type (every matrix is similar to an upper
triangular matrix is the equivalent formulation).  The same holds
for any connected algebraic group over an algebraically closed
field acting its on flag variety --- every element
is contained in a Borel subgroup and all
Borel subgroups are conjugate.

Derangements have proved to be very useful. In particular,
they have applications to images of rational points for
maps between curves over finite fields (and more generally
to higher dimensional varieties as well).  See \cite{gurwan}
for more details.  They also are useful in studying
probabilistic generation (see \cite{GLSS}).
Derangements are also useful in
giving bounds for the convergence rates of random walks on finite
groups \cite{Dia}.
We will explore these ideas further in future articles.

Recall that $G$ is called a Frobenius group of degree $n$ if $G$
acts transitively on a set of cardinality $n$ such that no
element in $G$ fixes $2$ points (and $|G| > n$). Surprisingly, it
was only very recently that Cameron and Cohen proved:

\begin{theorem}\label{1/n} \cite{cc}
If $X$ is a transitive $G$-set of cardinality
$n > 2$, then $\delta(G,X) \ge 1/n$ with equality if and
only if $G$ is a Frobenius
group of cardinality $n(n-1)$ (and in particular, $n$ is a prime power).
\end{theorem}

The proof is quite elementary.
Another proof of this is given in \cite{gurwan}
and the result is extended in various ways.  In particular, it was shown that:

\begin{theorem}\label{2/n} \cite{gurwan} If $X$ is a transitive
$G$-set of cardinality
$n > 6$, then $\delta(G,X) > 2/n$ unless  $G$ is a Frobenius
group of cardinality $n(n-1)$ or $n(n-1)/2$
(and in particular, $n$ is a prime power).
\end{theorem}

The proof of this result seems to require the classification of almost
simple $2$-transitive groups
(and so the classification of finite simple groups).

Note that when trying to produce lower bounds for
the proportion of derangements, there is no loss in assuming
that $G$ acts faithfully on $X$.  We will typically make that assumption.

We particularly want to focus on the case of primitive permutation groups
and simple and almost simple groups.  The Aschbacher-O'Nan-Scott theorem
\cite{aschscott}
gives the structure of primitive permutation groups and reduces many questions
about them to almost simple groups (groups which have a unique minimal normal
subgroup which is nonabelian simple).

A primitive permutation group $G$ of degree $n$ is
called affine if it preserves an affine
structure on the set.  This is equivalent to saying that $G$ has a nontrivial
normal elementary abelian $p$-subgroup $N$ for some prime $p$.  Necessarily,
$|N|=n$ is a power of $p$.  In particular, primitive Frobenius groups
are always
affine permutation groups.

The major part of this paper deals with $\delta(G,X)$ when $G$ is a finite
nonabelian simple group.  In particular, we will discuss a recent result
of the authors proving  a conjecture that has been attributed to Shalev.
This theorem is proved in  a series of papers by the authors starting with
this one -- see also \cite{FG1}, \cite{FG2} and \cite{FG3}.

\begin{theorem}\label{shalevconj}
There exists a positive number $\delta$ such that
$\delta(G,X) > \delta$ for all finite nonabelian simple groups
$G$ and all nontrivial transitive $G$-sets $X$.
\end{theorem}

Note that it suffices to prove the previous theorem when $X$ is a primitive
$G$-set (for if $f:Y \rightarrow X$ is a surjection of $G$-sets, then
$\delta(G,Y) \ge \delta(G,X)$). We also note that
as stated this is an asymptotic result -- we only need to
prove that there exists a $\delta > 0$ such that
for any sequence $G_i,X_i$ with $|X_i| \rightarrow \infty$
with $\delta(G_i,X_i) > \delta$ for all sufficiently large $i$.
This result is known for $G$ alternating
essentially by \cite{D} and \cite{lupyber} and for $G=PSL(d,q)$ for
a fixed $q$ \cite{shalev} for many families of actions.
Shalev's method used some difficult results about the order of a
random matrix; we use simpler properties of random matrices.

This paper is a partially
expository paper regarding variations on this theme
in \cite{FG1} and \cite{FG2}.  We will discuss in details analogous results
for algebraic groups and give the proof for finite Chevalley 
groups of bounded rank.

We will prove much more specific theorems and obtain much better
asymptotic results.  Our proof shows that we can take
$\delta$ to be roughly $1/25$ aside from finitely many exceptions
(and it is likely that there are no exceptions).

As in \cite{lupyber} and \cite{shalev}, we obtain
results about families of subgroups as well.  The following
result is proved in \cite{FG1}, \cite{FG2} and \cite{FG3}.

\begin{theorem}\label{classical}  Let $X:=X_n(q)$ be a classical 
group of dimension $n$
over $\F_q$.  Let $I(X)$ be the union of all proper
irreducible subgroups
of $X$ except when $X=Sp_{2m}(q)$ we do not include the irreducible
subgroups containing $\Omega^{\pm}_{2m}(q)$ with $q$ even.
Then $\lim_{n \rightarrow \infty} |I(X)|/|X| = 0$.
\end{theorem}

In the theorem, we can take $X$ to be a simple  classical group
or the full conformal subgroup or anything in between.  Moreover,
we can allow a center.
In \cite{FG4}, we will use this result to obtain some new results
about probabilistic generation.

One might think that Theorem \ref{shalevconj} is 
valid  for almost simple groups.
However, examples constructed in \cite{FGS} and \cite{GMS}
(coming from problems in coverings
of curves) show that the result fails for almost simple groups.
We give some examples of this phenomenon later in the paper.
The presence of field automorphisms is critical in producing
such examples.

However, we do prove that the result on simple groups does lead
to a weaker bound for primitive groups which are not affine.
See \S \ref{primitivesec}.

\begin{theorem} \label{primitive} Let $X$ be a primitive
$G$-set with $|X|=n$.
There exists a positive number $\delta$ such that
either
\begin{enumerate}
\item $\delta(G,X) > \delta/\log(n)$; or
\item $G$ preserves an affine structure on $X$.
\end{enumerate}
\end{theorem}

We will investigate the affine case in future work.

We will also consider a slight refinement of this problem:

Let $G$ be a normal subgroup of $A$ with $A$ and $G$ acting
transitively on a finite $A$-set $X$.  Assume that $A/G$ is
generated by the coset $aG$.

We wish to investigate the quantity
$\delta(A,G,X)$, the proportion of derangements in the coset $aG$.
We note the following easy facts:

\begin{enumerate}
\item the quantity $\delta(A,G,X)$ does not depend
on the choice of the generating coset $aG$; and
\item $\delta(G,X)=\delta(G,G,X)$.
\end{enumerate}

This quantity is important in studying
maps of varieties over finite fields via a Cebotarev density
theorem (see \cite{gurwan}
for more details).  In contrast to the case $A=G$,
there may be no derangements
in a given coset.  This turns out to be a very special and
important case in the study
of exceptional covers \cite{FGS} and graph theory  \cite{GLPS}.
See \S \ref{exceptional} and \S \ref{coset} for further discussion.

In \cite{gurwan}, the following was shown via an elementary proof:

\begin{theorem} \label{gwthm}Let $G$ be a normal 
subgroup of $A$ with $A/G$ cyclic.
Let $X$ be a transitive $A$-set of cardinality $n > 2$.
Either $\delta(A,G,X)=0$ or $\delta(A,G,X) \ge 1/n$. Moreover,
equality holds if and only if $A=G$ is a Frobenius group of order
$n(n-1)$.
\end{theorem}

One of the eventual goals of this project is to greatly improve this result.

We now give a brief sketch of the contents and ideas of the paper.

Let $G$ be a group acting primitively and faithfully on the transitive
$G$-set $X$ and let $H$ be the stabilizer of a point $x \in X$.
Set
$$\C_G(H):= \cup_{u \in G} uHu^{-1}.$$
An element $g \in G$ is  a derangement if and only if
$g \notin \C_G(H) $.

The proofs of many of these results are
heavily dependent upon the classification of finite
simple groups -- both in the fact that we are assuming the complete
list of finite simple groups and in using information about subgroups
of the finite simple groups.  Since Theorem \ref{shalevconj}
is really an asymptotic result,  
we are considering the following situation -- we have a sequence
$(G_i,X_i)$ where $G_i$ is a finite nonabelian simple group and
$X_i$ is a primitive $G_i$-set of cardinality of $n_i$.  We may assume
that $|G_i|$ (or $n_i$) tends to infinity.  We need to show
that $\liminf \delta(G_i,X_i) \ge \delta$ for some positive $\delta$
(a single $\delta$ for all such sequences).  This implies that
$\delta(G,X) \ge \delta$ for  all but finitely many simple  $G$
and primitive $X$, whence $\delta(G,X)$  is bounded away from $0$
for all simple  $G$ and primitive $X$.

In \cite{FG1}, \cite{FG2} and \cite{FG3}, we obtain much stronger results.

By passing to infinite subsequences, to prove Theorem \ref{shalevconj},
it suffices to
assume that all the $G_i$ are alternating groups (of increasing
degree), are
all Chevalley groups of a given type (and rank) over
fields of size $q_i$ with $q_i \rightarrow \infty$ or
are classical groups of dimension $d_i$ over fields of cardinality $q_i$
with $d_i \rightarrow \infty$.

In the case of alternating (and symmetric groups), we can apply the work
of Dixon \cite{D} and Luczak-Pyber \cite{lupyber}.  We improve some of
these results in \cite{FG1} and \cite{FG2}.
In the case of Chevalley groups of fixed type, we can use the theory of
algebraic groups and algebraic geometry to obtain the desired results.
Here the dichotomy is between subgroups that contain maximal tori and those
that do not. See \S \ref{algebraic}, \S \ref{boundedrank}
and \S \ref{maxranksec}.

Now consider the case that the $G_i$ are classical groups of dimension $d_i$
over a field of size $q_i$.  We subdivide this case further.  First of all,
either the $q_i \rightarrow \infty$ or we may assume that $q_i=q$ is constant.
In the first situation, by \cite{GL}, it suffices to consider only semisimple
regular elements.  We subdivide each case further using the idea of
Aschbacher's
classification of maximal subgroups of classical groups \cite{aschmax}.
In particular,
we consider subspace stabilizers and show that we can reduce
certain questions to
the study of the Weyl group (and so to symmetric groups).
We  prove in \cite{FG1} the following result about subspace stabilizers.
For the next theorem, in the case of a linear group, all subspaces
are considered degenerate.

\begin{theorem}  Let $G_i$ be a sequence of classical groups
with the natural module of dimension $d_i$.  Let
$X_i$ be a $G_i$-orbit of either totally singular or nondegenerate
subspaces (of the natural module) of dimension $k_i \le d_i/2$.
If $k_i \rightarrow \infty$, then $\lim \delta(G_i,X_i)=1$.
If $k_i$ is a bounded sequence, then there exists
$0 < \delta_1 <  \delta_2 < 1$
so that $\delta_1 <  \delta(G_i,X_i) < \delta_2$.
\end{theorem}

One of the key ingredients in this analysis is getting estimates for
the number
of conjugacy classes for finite Chevalley groups of rank $r$ over a
field of size
$q$.  We show that there is an explicit
universal constant $C$ so that the
number of conjugacy
classes of such a group is at most $Cq^r$
-- see \S \ref{conjugacyclasses}.  
See also Gluck \cite{gluck}
and Liebeck-Pyber \cite{liepyber} for weaker estimates.  These
results are of independent
interest and should be useful (see the above mentioned references
for some applications).
Two other important ideas in the proof are an upper bound for
the maximum size
of a conjugacy class and a result that says that
most elements in a classical group are nearly regular semisimple
(i.e. they are regular semisimple on
a subspace of small codimension).  Another ingredient
we use in the proof of Shalev's conjecture (for $q$
fixed) is precise estimates on proportions of
regular semisimple elements  (proved
via generating functions). See \cite{FNP}.
Finally, we require results on random permutations.

This article is organized as follows:

We first discuss derangements in algebraic groups (in algebraic actions)
-- see \S\ref{algebraic}.
We then prove Theorem \ref{shalevconj} for groups of bounded rank
in \S \S \ref{boundedrank}, \ref{maxranksec}  -- 
the second of which focuses on subgroups containing
a maximal torus.  As a corollary (\S\ref{generation}), we
solve for bounded rank  groups a problem studied by Dixon
\cite{D} and McKay (unpublished) for symmetric groups. The case of classical
groups with rank going to $\infty$ is treated in \cite{FG4}.
In \S\S \ref{exceptional}, \ref{coset}, we 
give some examples and mention the connection with so
called exceptional permutation actions and give a short proof
of Theorem \ref{gwthm}.  In \S\ref{primitivesec}, we 
then show how Theorem \ref{primitive}
follows as a corollary to Theorem \ref{shalevconj}.
In the final section, we tabulate some of our results about
conjugacy classes for classical groups.

\section{Algebraic Groups} \label{algebraic}

In this section, we investigate the existence of derangements in
(algebraic) permutation actions for connected algebraic groups.
We refer to \cite{alg-groups} for the basic results about algebraic
groups.

We first make a simple observation that holds for solvable
groups (not just algebraic groups).

\begin{lemma}  Let $G$ be a solvable group and $H$ a proper subgroup
of $G$.  Then $\cup_{g \in G} H^g \ne G$.
\end{lemma}

\begin{proof}  Let $A$ be the last term in the derived
series of $G$.  If $HA \ne G$, we can pass to $G/A$
and the result follows by induction on the derived length
(the case of abelian groups being obvious).  So assume that
$G=HA$ and in particular, $H$ does not contain $A$.
Then $H \cap A$ is normal in $G$ (since $H \cap A$ is normal
in $H$ and in the abelian group $A$).   It follows that
the only elements in $A$ which have fixed points are the
elements of $A \cap H$, a proper subgroup of $A$.
\end{proof}

We now consider connected algebraic groups.  We restrict
attention to semisimple groups.

Let $G$ be a connected semisimple algebraic group and $X$ a nontrivial
faithful algebraic $G$-set
$G/H$.  In particular, $H$ is a proper closed subgroup of $G$.
Let $J = \cup_{g \in G} H^g$.  Let $J'$ denote the complement of $J$ --
thus, $J'$ is precisely the set of derangements.

Let $T$ denote a maximal torus of $G$ and $N$ its normalizer.
Note that $T$ is self centralizing and $N/T$ is a finite group
(the Weyl group of $G$).  Moreover, any two elements of $T$
are conjugate in $G$ if and only if they are conjugate in $N$.
We recall that any two maximal tori of $G$ are conjugate.

\begin{lemma} \label{dense1} Let $G$ be a connected semisimple
algebraic group over an algebraically closed field.
Let $H$ be a closed subgroup of $G$ and $J= \cup_{g \in G} H^g$.
Then  the closure  ${\bar J}$ of $J$ is all of $G$ if and only if $T$
has a fixed point on $X$.
\end{lemma}

\begin{proof}  Since the set of semisimple elements of $G$
contains an open subvariety of $G$, the reverse implication is clear.

Assume that $J$ is dense in $G$.  Let $S$ be a maximal torus
of the connected component $H_0$ of $H$.  Let $d=|H:H_0|$.

Since $J$ is the image of
the morphism $f: H \times G \rightarrow G$ with $f(h,g)=h^g$,
it follows that $J$ contains a dense open subset of its
closure and so under this hypothesis an open subset of $G$.

Note that if $g \in G$ is semisimple regular, then the there
are at most $d^r$ solutions to $x^d=g$ (where $r$ is the rank of $G$)
-- for $x^d=g$ implies that $x \in C_G(g)$, a maximal torus of rank $r$.
This implies that the $d$th power map on $G$ is dominant and so
the set of $d$th powers of elements in $J$ also contains
an open subvariety of $G$.  This implies that the union of
the conjugates of $H_0$ contains an open subvariety
of $G$.   Since the union of the conjugates of
$S$ contains an open subvariety of $H_0$, we have that
$\cup_{g \in G} S^g$ contains an open subvariety of $G$.

By conjugating, we may assume that $S \le T$.
We have the surjection from 
$G/T \times S \rightarrow \cup_{g \in G} S^g$ given by
$(gT,s) \rightarrow s^g$, whence
$$ \dim G = \dim \cup_{g \in G} S^g \le \dim G + \dim S - \dim T, $$
and so $\dim S=\dim T$ and $S=T$.
\end{proof}

Of course, every element is contained in a Borel subgroup.  So if $H$
is a parabolic subgroup (i.e. an overgroup of a Borel subgroup), there
are no derangements.   We can give an easy proof that these are the only
examples if $H$ is connected.

\begin{theorem}  Let $G$ be a semisimple algebraic group over an algebraically
closed field $k$ of characteristic $p$.  Let $H$ be a closed proper subgroup
of $G$. Assume that $H$ is connected or that $p$ does not divide the order
of the Weyl group of $G$ (this includes the case $p=0$).
\begin{enumerate}
\item[(a)] If $H$ contains a maximal torus of $G$ and a 
regular unipotent element of
$G$, then $H$ is a parabolic subgroup of $G$.
\item[(b)] If $\cup_{g \in G} H^g =G$, then $H$ is a parabolic subgroup.
\end{enumerate}
\end{theorem}

\begin{proof}  If $\cup_{g \in G} H^g =G$, then previous lemma implies that
$H$ contains a maximal torus.  Clearly, it contains a regular
unipotent element, whence (b) follows from (a).

We now prove (a).  Let $H_0$ be the connected component
of $H$.  We will first show that $H_0$ contains a regular unipotent element.
We can then reduce to the case that $H$ is connected.

If $p=0$, this is clear because all unipotent subgroups are connected.
Let $T$ be a maximal torus of $H_0$ (which is also a maximal torus of $G$).
Since all maximal tori of $H$ are $H_0$-conjugate, it follows that
$H=N_H(T)H_0$, whence $|H:H_0|=|N_H(T):N_{H_0}(T)|$ is a divisor of the
order of the Weyl group of $G$.  In particular, it has order prime to $p$
and so $H_0$ contains all unipotent elements of $H$.

If $H_0$ is a parabolic subgroup, then so is $H$ (and indeed $H=H_0$
as any overgroup of a Borel subgroup is a parabolic subgroup).  So we
may assume that $H$ is connected.

Let $B_H$ be a Borel subgroup of $H$ containing $T$ and let $B$ be a Borel
subgroup of $G$ containing $B_H$.   Let $U$ be the unipotent radical of $B$.
Since $H$ is connected, $B_H$ contains a regular unipotent element $u$
as well (because every unipotent element of $H$ is conjugate to an element of
$U$).   We can write $u=v\prod U_{\alpha}(t_{\alpha})$ where the $\alpha$
are the simple roots relative to $T$, $t_{\alpha} \ne 0$ and $v \in [U,U]$.
It follows that $u^T[U,U]$ contains  all elements in $U$ which have
a nonzero entry in $U_{\alpha}$ for each simple root $\alpha$.  
Thus, $[u,T][U,U]=U$.  Since $U$ is nilpotent, this implies that
$U=[u,T]$ and so $B=TU \le H$. Thus, $B_H=B$
and $H \ge B$ as required.
\end{proof}

There are only a handful of examples of proper closed nonconnected
subgroups containing a conjugate of every element of $G$.
This requires
a result of Saxl and Seitz.  We note that the result of \cite{ss}
has the unneeded hypothesis that the characteristic is good (their proof never
uses this fact).  We will use the following fact in the next result --
any positive dimension closed subgroup of a simple algebraic group is
contained in a maximal closed subgroup.

\begin{theorem}  Let $G$ be a simple algebraic group over an algebraically
closed field $k$ of characteristic $p$.
Let $H$ be a closed proper subgroup of
$G$.  Assume that $H$ is not contained in a parabolic subgroup.
The following are equivalent:
\begin{enumerate}
\item[(a)] $H$ contains a maximal torus of $G$ and a
conjugate of every unipotent element of $G$;
\item[(b)] $H$ contains a conjugate of every element of $G$;
\item[(c)] The characteristic of $k$ is $2$ and
$(G,H)=(Sp(2m,k),O(2m,k))$ or $(G,H)=(G_2(k),A_2(k).2)$.
\end{enumerate}
\end{theorem}

\begin{proof}
Clearly (b) implies (a).  

We next show that (a) implies (c).  So assume that $H$ satisfies
(a). If $H$ is maximal (among closed
subgroups), then Theorem C of \cite{ss} shows that (c) holds.

Let  $M$ be a maximal closed subgroup of $G$
containing $H$. Then $(G,M)$ satisfies the conclusion of (a) as well
and so as noted, $(G,M)$ satisfies (c).  In particular, $k$
has characteristic $2$. Moreover, $H$ must have
maximal rank and is not connected.  

If $G=G_2(k)$ and $H$ is a proper (disconnected) rank $2$ subgroup
of $M$, then the only possibility is that $H$ is contained in the
normalizer of a maximal torus, which does not contain a conjugate
of every unipotent element.

So we may assume that $G=Sp(2m,k)$.  If $m=1$, then $M$ is the normalizer
of a maximal torus $T$ and $M/T$ has order $2$ and so clearly $H=M$.
So consider the case that $G=Sp(2m,k), m \ge 2$.  If $H$ acts reducibly on
the natural module $V$ for $G$, then $H$ is contained in the stabilizer
of a proper subspace $W$.  Take $W$ of minimal dimension.
Since $H$ is not contained in a parabolic subgroup, it
follows that this subspace is nondegenerate.
The stabilizer of such a subspace is not contained
in a conjugate of $M$, a contradiction.

Suppose that the connected
component $H_0$ does not act irreducibly on $V$.  Then either $H$ is contained
in a maximal subgroup not contained in $O(2m,k)$, a contradiction or
$V=W_1 \oplus W_2$, where $H_0$ is irreducible on each $W_i$ and $W_i$ is a
maximal totally
singular subspace of dimension $m$.  Thus,
$H$ is contained in the stabilizer of a pair
of complementary totally singular subspaces.
We claim that $H$ contains no transvections.
A transvection in $H$ cannot swap the two spaces and
so would have to stabilize
each $W_i$.  The action on $W_1$ is dual to that on
$W_2$ and so the element is
not a transvection.

Thus, $H_0$ acts irreducibly on $V$, whence $H_0$ is semisimple.
Since $H$ has rank $m$, this forces $H_0$ to contain the connected
component of $O(2m,k)$, whence $H=O(2m,k)$.

All that remains is to verify that (c) implies (b).

This is well known for the first family (see \cite{ss}, Lemma 4.1)
The latter case is an easy consequence of the first case 
(since $G_2(k) \le Sp(6,k)$,
$Sp(6,k)=O(6,k)A_2.2$ and $A_2.2 = G_2(k) \cap O_6(k)$).
\end{proof}

\section{Groups of Bounded Rank I}\label{boundedrank}

In this section and the next, we consider the case where the 
groups have bounded rank.
We will prove Theorem \ref{shalevconj} in this case. 
The next section deals with subgroups containing a maximal torus.
We deal with the other cases in this  section.

 As we have
observed, as stated it is an asymptotic result.  We only need to produce
a $\delta$ so that the proportion of derangements is at least $\delta$
for all  but finitely many cases. If this fails, there would be a sequence
with the proportion of derangements all less than $\delta$. Thus,
Theorem 1.3 is an asymptotic result (as noted, eventually we want a 
non-asymptotic
version).  Since the groups have bounded rank, we may assume that they have
fixed type $X(q_i)$ with $X$ a simple algebraic group and only the 
field size is varying.
We can use methods of algebraic geometry and algebraic groups to study this
situation.

We recall that  $F^*(H)$ is the generalized Fitting subgroup of $H$.
See \cite{aschbook}.  In particular, $F^*(H)$ simple just means that
$F^*(H) \le H \le \aut(H)$.  There is no harm in considering covering
groups of almost simple groups since all the maximal subgroups will
contain the center.

Fix a type of simple algebraic group $X$ of rank $r$.  Let
$\sigma$ be an endomorphism of $X$
with fixed point group $X_{\sigma}$ of finite order.  We will typically
write $X(q)$ if $q$ is the absolute value of the eigenvalues of $\sigma$
on the character group of the maximal torus $T$ of $X$. In  the case of
the Suzuki or Ree groups $q$ will not be an integer.  This
will cause no problems.  Indeed, in those cases, one knows all the maximal
subgroups and it is quite easy to obtain our results.  We may take $X$
simply connected or of adjoint type or anything in between -- this allows us
to  obtain results for Chevalley groups generated by inner-diagonal 
automorphisms.

The maximal subgroups $H$ of $X(q)$ are of four types:

\begin{enumerate}
\item $|H| < N$ for some fixed $N=N(X)$;
\item $H = N_{X(q)}(X(q'))$ for some $q'$ dividing $q$
(this includes the twisted forms, e.g., ${^2}E_6(q) \le E_6(q^2)$);
\item $H= Y_{\sigma}$ where $Y$ is a proper $\sigma$-invariant semisimple
algebraic subgroup of $X$ of rank
$s < r$;
\item $H= Y_{\sigma}$ where $Y$ is a proper $\sigma$-invariant
algebraic subgroup of $X$ of maximal rank $r$.
\end{enumerate}

This is well known for the case of classical groups (see \cite{aschmax}).
If $H$ is of exceptional type, this follows in a very precise way
from results of Liebeck and Seitz \cite{LS1}.   See
the remarkable paper of Larsen and Pink \cite{larsenpink} for a classification
free proof of the previous result.  Note that if $Y$ is a maximal
$\sigma$-invariant positive dimensional algebraic subgroup, then
either $Y$ has maximal rank or the connected component of $Y$ is semisimple 
(for if the unipotent radical of $Y$ is nontrivial, $Y$ is a parabolic subgroup
by the Borel-Tits theorem and if it contains a normal torus, then
$Y$ contains a maximal torus).

We also note the following well known result (see also \cite{LS1}).

\begin{lemma}  If $X$ is a simple connected algebraic group and
$Y$ is a proper positive dimensional $\sigma$-invariant subgroup , 
then $Y$ is contained in a maximal proper closed $\sigma$-invariant subgroup.
Moreover, there is a bound $m=m(X)$ for the number of connected components
for any maximal $\sigma$-invariant  closed subgroup.
\end{lemma}

Let ${\mathcal M}_i$ denote the maximal subgroups in the 
corresponding families
$i=1,2,3$ or $4$ above.

We will deal with each of these families separately.  In this section, we deal
with the first three cases.  In the next section, we deal with 
the remaining case.
The purpose of this section is to prove:

\begin{theorem} \label{m123}  $\lim_{ q \rightarrow \infty} 
|\cup_{i =1}^3 \cup_{M \in {\mathcal M}_i} M|/|X(q)| = 0$.
\end{theorem}

This is not true for ${\mathcal M}_4$ (compare with the results 
for algebraic groups).
We first start with some general known results.

\begin{lemma} Let $x \in X(q)$. Let $C$ be the centralizer of
$x$ in $X(q)$.
\begin{enumerate}
\item If $x$ is unipotent, 
then $|C|$ is divisible by $q^r$.
\item $|C| \ge (q-1)^r$.
\end{enumerate}
\end{lemma}

\begin{proof}  If $x$ is unipotent, the result follows since
all unipotent classes are known as well as their centralizers.
One  can also give a proof based on fixed points on unipotent
groups.

We prove the second  statement more generally for reductive
groups of rank $r$.
Write $x=su=us$ with $u$ unipotent and $s$ semisimple.
Pass to the connected component of $D$ of $C_X(s)$.  This is still reductive
of rank $r$.  Write $D=AB$ with $A$ a central torus in $C$ and
$B=[C,C]$ semisimple with $A$ of rank $a$ and $B$ of  rank $b$.
Since  a torus of rank  $a$ over the field of $q$ elements has
at least size $q-1$ and $C_{B(q)}(u)$ has order divisible by $q^b$, we
see that $|C| \ge (q-1)^aq^b$, whence the second statement holds.
\end{proof}

The following was originally proved by Steinberg
in the case of simply connected $X$.  See \cite{Ca} or  \cite{ssclasses}.

\begin{lemma}  The number of conjugacy classes of semisimple
elements in $X(q)$ is at most $q^r$ with equality if $X$ is
simply connected.
\end{lemma}

The next result follows from \cite{GL}.

\begin{lemma}  The proportion of regular semisimple elements in $X(q)$
is greater than $1 - 5/(q-1)$.
\end{lemma}

The previous result indicates that the proportion of elements which are not
semisimple regular goes to $0$ linearly with $1/q$.  The same is thus true
for the  set of derangements which are not semisimple regular.
Thus, it suffices to consider the set of derangements which are semisimple
(and indeed regular).  We will do so in the next two sections without further
comment.

\begin{lemma}  $|\cup_{M \in {\mathcal M}_1} M |/|X(q)| \rightarrow 0$
as $q \rightarrow \infty$.
\end{lemma}

\begin{proof} $ \cup_{M \in {\mathcal M}_1} M$ is a union of at most
$N'$  conjugacy classes for some $N'$ (that depends only
on $N$ and so  only on $X$). Thus  the union has
order at most $|X(q)|N'/(q-1)^r$ and the result follows.
\end{proof}

\begin{lemma}  $|\cup_{M \in {\mathcal M}_2} M |/|X(q)| \rightarrow 0$
as $q \rightarrow \infty$.
\end{lemma}

\begin{proof}  Consider $X(q')$. The number of semisimple conjugacy
classes in $X(q')$ is at most $(q')^r$.
Let $S(q',q)$ denote the union of the semisimple conjugacy classes
of $X(q)$ intersecting $X(q')$. Thus,
$$
|S(q',q)|\le |X(q)|(q')^r/(q-1)^r.
$$

In the case of the Suzuki or Ree groups, we write
$X=X(p^{2a+1})$ (this conflicts slightly with our notation
above).  The number of possible classes of subfield groups
is the number of distinct prime divisors of $2a+1$, whence
the estimate above shows that the union of the semisimple
elements in any subfield group is certainly at most
$\sum_b |X(q)|q^{r/b}/(q-1)^r$, where $b$ ranges over
prime divisors of $2a+1$.  This yields the result.

Consider the remaining cases.
Write $q=p^a$.  Note that for each choice of $q'$ (corresponding
essentially to a prime divisor of $a$), there are at most $2c$
choices for $S(q,q')$ where $c$ is the order of the group
of outer diagonal automorphisms  ($6c$ in case $X=D_4$).
This is because
we may take $\sigma = \alpha \tau f_{q'}$ where
$\tau$ is a graph automorphism,
$\alpha$ is a diagonal automorphism
and $f_q$ is the standard Frobenius (any two such elements in
the coset with the same order are conjugate up to diagonal
automorphisms -- see I.7.2 \cite{gorly}).  In fact as noted above,
we only need to consider
semisimple elements, the diagonal outer automorphisms will not make
a difference, but we do not need to use this.
\end{proof}

\begin{lemma}  $|\cup_{M \in {\mathcal M}_3} M|/|X(q)| \rightarrow 0$
as $q \rightarrow \infty$.
\end{lemma}

\begin{proof}  It follows by the theory of high weights
if $X$ is classical \cite{GKS} and by \cite{LS2} if $X$ is exceptional
that there are only
finitely many conjugacy classes (with a bound depending only upon $X$)
in ${\mathcal M}_3$.  Thus, it suffices to show the result for
a fixed type of subgroup $Y < X$.  Then $Y_{\sigma}$ has at most
$cq^s$ conjugacy classes of semisimple elements (where $c$ is the
number of connected components of $Y$ -- note that $c$ is bounded
in terms of  $X$).  It follows
that $|\cup_{g \in X_{\sigma}} Y_{\sigma}^g| \le |X_{\sigma}|cq^s/(q-1)^r$,
whence the result.
\end{proof}

This completes the proof of Theorem \ref{m123}.  The next section
deals with ${\mathcal M}_4$.

\section{Maximal Rank Subgroups}\label{maxranksec}

In this section we consider ${\mathcal M}_4$.  It follows from the results
on algebraic groups that the proportion of derangements will be
positive in this case.  For the Suzuki and Ree groups, one just inspects
the maximal rank subgroups and the result about derangements follows quite
easily.  We assume for the rest of the section that we are not 
in any of those cases.  We remark again that it suffices
to consider only regular semisimple elements (since as $q \rightarrow
\infty$, the proportion of regular semisimple elements is
$1 + O(1/q)$).

Keep notation as in the previous section.  Let $Y$ be
a $\sigma$-stable subgroup of $X$ of maximal rank
and $H=Y_{\sigma}$.
The possibilities are that $Y$ is a parabolic subgroup
(maximal with respect to being $\sigma$-stable) or
is reductive.  Let $Y_0$ denote the connected component of $Y$.
Let $H_0=(Y_0)_{\sigma}$.  This is a normal subgroup of $H$.

There exists
a $\sigma$-stable maximal torus $T$ contained in a Borel
subgroup $B$ of $X$.
A maximal torus of $X_{\sigma}$ is $S_{\sigma}$ where $S$
is a $\sigma$-stable maximal torus of $X$.  There is a notion
of nondegenerate maximal tori (for example, if $X=SL(n)$, then
over the field of $2$-elements, a maximal torus might be trivial,
see \S3.6 in \cite{Ca} for details).
We will just note that if the maximal torus contains a regular
semisimple element, then
$N_{X_{\sigma}}(S_{\sigma})=N_X(S)_{\sigma}$ -- this follows
since $S = C_X(S_{\sigma})$.   Moreover (for fixed $X$), if
$q$ is sufficiently large, all maximal tori contain regular
semisimple elements (indeed almost all elements are regular
semisimple).

Let $W$ be the Weyl group of $G$ (more precisely identify
$W = N_X(T)/T$).  Consider the semidirect product
$W\langle \sigma \rangle$.  There is a bijection between
conjugacy classes of maximal tori in $X_{\sigma}$ and
$W$-classes of elements in the coset $\sigma W$
(see \cite{ss2} or \cite{Ca}). In particular, if $\sigma$ is a field
automorphism, $\sigma$ commutes with $W$ and so the
correspondence is with $W$-conjugacy classes (this latter fact is still
true for all groups of type $A$ and type $D_n$ with $n$ odd).

Let $T_w$ denote a maximal torus of $X_{\sigma}$ corresponding
to $\sigma w$.  Let $N_w$ be the normalizer in $X_{\sigma}$ of
$T_w$.  Then $|N_w:T_w|=|C_W(\sigma w)|$.
Let  $f(w)$ be the size of $W$-class of $\sigma w$.
So $f(w)=|W:C_W(\sigma w)|=|W||T_w|/|N_w|$.

In particular, we see that
$$
|\cup_{g \in X_{\sigma}} T_w^g|/{|X_{\sigma}}| < |T_w|/|N_w|
= f(w)/|W|.
$$

Since a semisimple regular element lies in a unique maximal
torus, it follows that the union of all regular semisimple
elements of $X_{\sigma}$ that are conjugate to an element
of $T_w$ has cardinality at most $|X_{\sigma}|f(w)/|W|$.

Since the proportion of elements which are not semisimple
regular tends to $0$ as $q \rightarrow \infty$ and the same is true
for each maximal torus, it follows that the inequality
above becomes equality as $q \rightarrow \infty$.

 We first show
that the collection of elements which are conjugate to
an element of $H$ but not $H_0$ is small.
We need the following result.  A very easy result
gives an upper bound (always at least $1/2$) for the proportion of derangements
contained in $H_0$ (assuming that $H \ne H_0$).

\begin{lemma}  Let $G$ be a connected reductive algebraic group
with $\sigma$ an endomorphism of $G$ such that
$G_{\sigma}$ is finite.  Assume that all eigenvalues
of $\sigma$ on the character group of $T$ have
absolute value $q$. Let $S$ and $T$ be distinct
$\sigma$-stable maximal tori of $G$.  Then
$|T_{\sigma}:(S \cap T)_{\sigma}| \ge (q-1)/2$.
\end{lemma}

\begin{proof} Consider a counterexample with $\dim G$ minimal.
Since $G$ is reductive,
$G$ is the central product of $Z$ and $H$ where $H$ is semisimple
and $Z$ is the connected component of the center of $G$.  Since
$Z$ is contained in every maximal torus, there is no loss
in taking $G=H$ to be semisimple.   We can replace $G$ by its
universal central extension (since the center will be contained
in every maximal torus) and so assume that $G$ is a direct
product of simply connected simple algebraic groups.

If $S \cap T = Z(G)$, the result is clear (pass to the simple
case).  Otherwise, we can consider
$H=C_G(x)$ with $x \in S \cap T \setminus{Z(G)}$.
Then $H$ is connected and reductive and $S,T$ are maximal tori in $H$.
\end{proof}

Note that if $G=SL(2)$, we do have equality in the previous result.

\begin{prop}\label{disconnected} 
$$
\lim_{q \rightarrow \infty}
|\cup_{g \in X_{\sigma}}(Y \setminus Y_0)_{\sigma}^g|/|X_{\sigma}|
= 0.
$$
\end{prop}

\begin{proof}  It suffices to consider a single coset $yY_0$ for
some element $y \in Y_{\sigma} \setminus Y_0$.

We will obtain an upper bound on the number of conjugacy classes
of semisimple regular elements of $X_{\sigma}$ that intersect
$yY_0$.  We will do this by bounding the number of $Y_{\sigma}$
classes in that coset.

Suppose that $u \in yY_0 \cap X_{\sigma}$ is a
semisimple regular element. Then the centralizer
of $u$ in the algebraic group is a $\sigma$-stable
maximal torus $T$.  Let $S$ be a $\sigma$ stable
maximal torus of $Y_0$ containing $T \cap Y_0$.
the number of $(Y_0)_{\sigma}$ conjugates of $u$
is
$$
|(Y_0)_{\sigma}:(S \cap T)_{\sigma}| \ge |(Y_0)_{\sigma}|(q-1)/2|S_{\sigma}|.
$$

Since $|S_{\sigma}| \le (q+1)^r$, it follows that
the  number of conjugates of $u$ in the coset $u(Y_0)_{\sigma}$
is at least  $|(Y_0)_{\sigma}|q^{r-1}/2$ (up to a small error term).
This implies that there are at most $2q^{r-1}$ classes of semisimple
regular elements in this coset (again up to a term of smaller order).
Since each class has size approximately $O(|X_{\sigma}|)/q^r$,
the union of these class has size  $O(|X_{\sigma}|/q)$
as required.
\end{proof}

We now consider the connected component $Y_0$ and its fixed
points $H_0$.  We first note that if $H_0 \ne H$, then
we have the following easy estimate for derangements.

\begin{lemma}  If $H \ne H_0$, then 
$|\cup_{g \in X_{\sigma}} H_0^g|/|X_{\sigma}|< 1/|H:H_0|
\le 1/2$.
\end{lemma}

\begin{proof}  Since $H$ normalizes $H_0$,
$\cup_{g \in X_{\sigma}} H_0^g$ is 
the union where $g$ ranges over a transversal
of $X_{\sigma}/H$, whence the cardinality of this union is less than
$|X_{\sigma}:H||H_0| =1/|H:H_0|$.  
\end{proof}

We just remark that Lang's theorem implies that $|H:H_0|$
is the number of $\sigma$-stable cosets of $Y_0$ in $Y$.

Let $S$ be a $\sigma$ stable maximal torus
of $Y_0$.  Then $S=xTx^{-1}$ where $x^{-1}\sigma(x) \in N(T)$.
Note that $x^{-1}N(S)x = N(T)$.  So we have subgroups
$T \le x^{-1}N_H(S)x \le x^{-1}N_Y(S)x \le N(T)$ and this
gives rise to corresponding subgroups $1 \le W_0 \le W_1 \le W$
in $W$ the Weyl group of $T$.

The $\sigma$-stable maximal tori of $H$ (up to $H_{\sigma}$-conjugacy) 
are of
the form $ySy^{-1}=yxT(yx)^{-1}$ where $v:=y^{-1}\sigma(y) \in N_H(S)$.
Moreover, we see that $S$ is conjugate to $T_w$ where
$$
w = (yx)^{-1}\sigma(yx)T = x^{-1}y^{-1}\sigma(y)\sigma(x)
\in \tau W_0,
$$
where $\tau =x^{-1}\sigma(x)T \in W$.

Thus, setting $R$ to be the union of all  $X_{\sigma}$ conjugates
of maximal tori of $H_{\sigma}$, we see that
$|R|/|X_{\sigma}|  \le \sum f(w)/|W|$
where the sum is a set of representatives $\Gamma$
of conjugacy classes in $W$ that are represented by elements
in $\tau W_0$.

We note that $R$ is precisely the set of conjugacy classes
of regular semisimple elements conjugate to an element of $H_{\sigma}$
(since the centralizer of such an element will be a maximal torus
in $H$).  Thus, we have an upper bound for the proportion
of regular semisimple
elements in $X_{\sigma}$ that are not derangements in
$X_{\sigma}/Y_{\sigma}$.  By Proposition \ref{disconnected}, 
we can replace $H_{\sigma}$ by $Y_{\sigma}$
(up to an $O(1/q)$ term) and we can consider all elements
(not just regular semisimple elements) by introducing
another such term (by \cite{GL}) and so we see that
$\delta(X_{\sigma},Y_{\sigma}) \ge 1 - \sum_{w \in \Gamma} f(w)/|W| + O(1/q)$.
We just need to bound $\sum_{w \in \Gamma} f(w)/|W|$ away from $1$.

There is one very easy case -- if $\sigma$ does not involve
a graph automorphism and $W_1 \ne W$, then
$1-\sum_{w \in \Gamma} f(w)/|W| \le \delta(W,W/W_1)$. Note
in this case $f(w)$ is just the size of the $W$-conjugacy class
of $w$.  This can
be computed for the exceptional Weyl groups.  In any case,
for bounded rank, we can even use the Jordan bound to see
this is bounded away from $1$.  For classical groups, using
\cite{D}, \cite{lupyber}, \cite{FG1}, \cite{FG2} we see that this
quantity will typically be at most $2/3$.

If $\sigma$ does involve a graph automorphism, then
we consider the group $Z$ defined above. 
Since $\sigma$ stabilizes both $W_1$ and $W_0$,
we can define $Z_1$ and $Z_0$ in an obvious manner.
Note that in this case $f(w)$ is the size
of the $W$-class of $\sigma w$.  Thus,
we still have a bound $1-\sum_{w \in \Gamma} f(w)/|W| \le
\delta(Z,W,Z/Z_1)$.

We note by inspection that unless $W=W_1$, there are
always derangements in the coset $\sigma W$.  Since
except for type $D_4$, $W$ can be thought as of a subgroup
of index $2$ in $Z$ (only the graph automorphism
makes a difference), exceptionality would force $|W:W_1|$
to have odd index.

There are only a few cases where $W=W_1$. 
In all cases, we see that whenever this happens
$H \ne H_0$ and so the upper bound of $1/2 + O(1/q)$ holds.
 One possibility is
that $H$ is a maximal torus.  In that case, we note directly that
$f(w)/|W| \le 1/2$.  Another possibility is $X=G_2$ and $H=A_2$.
Similarly, there is the possibility of $(X,Y)=(F_4,D_4)$ (in characteristic
$2$).  One sees that in these cases not all maximal tori are represented
in a given $H$.

The only other such possibility for $X$ classical is in characteristic $2$
 with $X$ of type $B_n$ and $Y$ of type $D_n$.  Then
$W_0$ has index $2$ in $W=W_1$.  In this case, one sees that
there are two possible forms of $Y_{\sigma}$ (i.e. two conjugacy
classes corresponding the single $X_{\sigma}$ class of $\sigma$
stable conjugates of $Y$) -- the two forms of orthogonal groups.
One form of the orthogonal group has maximal tori $T_w$ with
$w \in W_0$ and the other the complement (note a maximal torus
is contained in a unique orthogonal group in the symplectic group).
Thus, $\delta(Sp(2m,2^a),O^{\epsilon}(2m,2^a)) = 1/2 + O(1/2^a)$.

We note that our analysis works for any form of the Chevalley
group and for any fixed coset in the group of inner-diagonal
automorphisms.

Thus, we have proved:

\begin{theorem}\label{maxrank} Let $r$ be a positive integer.
Let $S$ be a simple Chevalley group of rank
at most $r$ over the field of $q$ elements
and $S \le G$ with $G$ contained in
the group of inner-diagonal automorphisms of $G$.
Let $X$ be a transitive faithful $G$-set.
Then there exists $\delta > 0$ such that
$\delta(G,S,V) \ge \delta + O(1/q)$ for any transitive
$G$-set $V$.
\end{theorem}

We will give an explicit $\delta$ in the sequel.
Note that the error term depends only on $r$ (and 
we do remove that dependence in \cite{FG1}, \cite{FG2} and \cite{FG3}).

If $X$ is classical, then the possibilities for $Y$ are
rather limited.  There is the special case in characteristic
$2$ when $X$ is symplectic and  $Y$ is an orthogonal
group.  The remaining cases are essentially when
$Y_{\sigma}$ is the group preserving a decomposition of
the space or $Y_{\sigma}$ is an extension field group
(both forms of $Y \le C \wr S_m$ where $C$ is a classical group
on a subspace).   In the bounded rank case, we have seen
that we could work with the connected piece.

If the rank increases, there are two added complications.
If $q$ is fixed, then we can no longer deal with only
semisimple regular elements (the error term may be larger
than the main term).  Even if $q$ increases, the error term
associated with reducing to the connected component
may be increasing with the rank.  Thus, the analysis is much
more difficult. See \cite{FG1} and \cite{FG2}.  
We also want to produce an explicit $\delta$
that is valid for either all cases or all but a specified
finite set of cases.

We close this section by considering a few examples.

\begin{enumerate}

\item Let $G=PSL(n,q)$ and let $H$ be the stabilizer of
a $k$-dimensional vector space.  In this case $H$ is the set of
fixed points of a connected subgroup
and so we see from the analysis above that for a fixed $n$,
$\lim_{q \rightarrow \infty} \delta(G,H) = \delta(S_n,Y_k)$
where $Y_k = S_k \times S_{n -k}$ is a Young subgroup.
By \cite{D}, $\delta(S_n,Y_k) \ge 1/3$ and by \cite{lupyber}
(and also by \cite{FG1}),  $\delta(S_n,Y_k) \rightarrow 1$
as $k \rightarrow \infty$ (for $k \le n/2$).  This example
holds more generally for any parabolic subgroup -- the limiting
proportion  of derangements is precisely the proportion of derangements
of the Weyl group acting on the 
cosets of the corresponding parabolic subgroup.

\item Let $G=Sp(2n,q)$ with $q$ even.  Let $H=O^{\epsilon}(2n,q)$.
Then $H$ is the set of $\sigma$ fixed points on some $\sigma$
invariant conjugate of $O(2n) < Sp(2n)$.  The Weyl group of
the connected component of $O(2n)$ has index $2$ in the
Weyl group of $Sp(2n)$ and so we see that
$\lim_{q \rightarrow \infty} \delta(G,H) = 1/2$
(each type corresponds to maximal tori in one coset
of the Weyl group of $\Omega$).
Note that a regular semisimple element is contained in
precisely one orthogonal group.

\item Let $G(q)=E_8(q)$ and $H(q)=D_8(q)$. 
Since the corresponding algebraic subgroup
is connected, it follows that
$$
\lim_{q \rightarrow \infty} \delta
(G(q),G(q)/H(q)) = \delta(W(E_8), W(D_8)).
$$

\end{enumerate}

\section{Generation and Derangements}\label{generation}

In this section, we indicate how some generation results follow immediately
from our results.  See \cite{FG4} for more results about
probabilistic generation that follow from the results in this paper,
\cite{FG1}, \cite{FG2} and \cite{FG3}.  

If $G$ is a finite simple group, set $P_G(x)$ the probability
that a random $y \in G$ satisfies $G=\langle x, y \rangle$.  Let $P_G$ be
the minimum of $P_G(x)$ over all nontrivial $x$.  It follows by \cite{GK} that
$P_G > 0$.  One can easily deduce the 
following result (a special case of \cite{GLSS}):

\begin{theorem}  Let $X$ be a type of simple algebraic group.
Then one has that $\lim_{q \rightarrow \infty} P_{X(q)}=1$.
\end{theorem}

 Recall that a group $G$ is generated {\it invariably} by the elements
$x_1,\cdots,x_m$ if the elements $y_1,\cdots,y_m$ generate $G$
whenever $y_i$ is conjugate to $x_i$ for $i=1,\cdots,m$. Luczak and
Pyber \cite{lupyber} proved the following conjecture of McKay, useful
in computational Galois theory.  (We make the constants in Theorem 5.2
more explicit in forthcoming work).

\begin{theorem} (\cite{lupyber}) There exists $N$ so that for all $n \geq N$
and all $\epsilon>0$, there is a constant $C=C(\epsilon)$ so that $C$
permutations, chosen from $S_n$ uniformly and independently, generate
$S_n$ invariably with probability at least $1-\epsilon$. \end{theorem}

For Chevalley groups of bounded rank, we have the following
result.

\begin{theorem} Let $X$ be a type of simple algebraic group.
For any $\epsilon>0$, there is a constant $C=C(\epsilon)$
(not depending upon $q$) so that $C$
elements, chosen from $X(q)$ uniformly and independently, generate
$X(q)$ invariably with probability at least $1-\epsilon$. \end{theorem}

\begin{proof}  
The probability that some $y_1, \ldots, y_m$ generate
a maximal subgroup in ${\mathcal M}_i, i \le 3$ tends to $0$ 
as $q \rightarrow \infty$ by Theorem \ref{m123}.
Indeed our proof shows that this probability is $O(1/q^m)$.
There are at most $d$ (depending only on $X$)  conjugacy classes
of maximal subgroups in ${\mathcal M}_4$ and the probability that some
conjugate of a random element $x \in G$ is contained in one is at
at most $1 - \delta$ for some $\delta > 0$ (for some $\delta$ depending
only on $X$). Thus, the probability that some collection of $y_i$
are contained in one of these maximal subgroups is
at most $d(1 - \delta)^m$.  So for $q$ sufficiently large,
we can choose an $m$ so with  probability greater than 
 $1-\epsilon$, $m$ random elements invariably generate $X(q)$.
If necessary, enlarge $m$ so that this is true for the finitely
many possibilities left over.  Note that we are ignoring the possibility
that $X(q)$ may not be simple -- however, this quotient is bounded in
terms of $X$ and so is not a problem.
\end{proof}

The analogous results for unbounded rank will be done
in a future article.

\section{Exceptionality and Derangements}\label{exceptional}

In this section, we discuss the notion of exceptional permutation
representations and its connection to curves.  See \cite{GMS}
for a more elaborate discussion of these ideas.

Let $G$ be a normal subgroup of $A$. Let $X$ be a transitive
$A$-set that is also transitive for $G$.
We say that
$(A,G,X)$ is exceptional if $A$ and $G$ have no nontrivial
common orbits on $X \times X$ (the trivial orbit being the
diagonal).
We note the following easy example.  See \cite{GMS} for more
examples and some classification theorems.

Recall that a Hall subgroup $H$ of a finite group $G$ is
a subgroup with $\gcd(|H|,|G:H|)=1$.

\begin{theorem} Let $A$ be a finite group and
$G$ a normal Hall subgroup. Then
$A=GD$ for some complement $D$ (by the Schur-Zassenhaus
theorem).  Let $H=N_A(D)$ and $X=A/H$.  Then
$(A,G,X)$ is exceptional.
\end{theorem}

\begin{proof}  Suppose not. We can identify $X$ with
the set of conjugates of $D$.  Let $K=C_G(D)=G \cap H$.
It is easy to see that exceptionality is equivalent
to $K$ and $H$ having no common orbit (other than
$D$ itself).  Suppose $D \ne E$ is in a common orbit.
Then the length of this orbit divides $|K|$ and
in particular  has order
prime to $|D|$.  Thus, $|H:N_H(E)|$ has order prime to
$D$, whence $N_H(E)$ contains a Hall $\pi$-subgroup $D_1$
of $H$ (where $\pi$ is the set of primes dividing $|D|$).
Since $D$ and $D_1$ are both Hall $\pi$-subgroups of $H$,
they are conjugate in $H$ (by the Schur-Zassenhaus theorem),
whence $D$ normalizes some $K$-conjugate of $E$.  So
we may assume that $D$ normalizes $E$.  Then $DE$ is a
$\pi$-subgroup, whence $D=E$, a contradiction.
\end{proof}

In particular, this result applies to the case
$G$ is a Chevalley group defined over the field of
$q^b$ elements, $b$ is relatively prime to  $|G|$
and we take
$D$ to be the cyclic group of order $b$ of field
automorphisms of $G$.  Specifically, take
$G=L(2,p^b)$ with $b$ any prime not dividing
$p(p^2-1)$.  See \cite{GMS}.  In this family 
of examples, the degree of the permutation representation
has size less than $p^{3b}$ and the proportion of derangements
is less than $1/b$ (since all derangements are contained
in $G$ and $|A:G|=b$).   Thus, even for almost simple groups acting
primitively, 
the proportion of derangements can be less than any given $\epsilon > 0$.
The best general result one could hope for is $\delta(A,X) > C/\log |X|$.
See \S \ref{primitivesec} for such a result.

Some easy facts about exceptional triples are
(see \cite{FGS}, \cite{gurwan}, \cite{GMS}):

\begin{enumerate}

\item If $A/G$ is generated by the coset $aG$,
then $(A,G,X)$ is exceptional if and only if every element
in the coset $aG$ has a unique fixed point or
equivalently $\delta(A,G,X)=0$.

\item If $A/G$ is cyclic and $(A,G,X)$ is exceptional, then so
is $(A,G,Y)$ where $Y$ is the image of a morphism
of $A$-sets from $X$ is an $A$-morphism.

\end{enumerate}

In particular, if $A/G$ has prime order and  $(A,G,X)$
is exceptional, then $\delta(A,X) < |G|/|A|$ (since all
derangements are contained in $G$).

So the analog of Shalev's conjecture for almost simple
groups fails.  In future work, we hope to obtain a result
that says that Shalev's conjecture holds except for certain
primitive actions (mostly related to the case where the
point stabilizer is the set of fixed points of some
Lang-Steinberg endomorphism of an algebraic group).

As we have remarked, $\delta(A,G,X)$ is related to
images of rational points for maps between curves and
higher dimensional varieties over finite fields.  The
connection is through the following estimate that
follows from the Cebotarev density theorem (see \cite{gtz}
or \cite{gurwan}for more details).

We make this more precise.

Let $U,V$ be smooth projective curves defined over $F:=F_q$
the field of $q$ elements.   Let $U(q^a)$ denote the $F_{q^a}$
rational points of $U$.
Let $f:U \rightarrow V$ be a separable rational map of
degree $n$ also defined over $F$.  Let $F(U)$ and $F(V)$
be the function fields of $U$ and $V$ over $F$.
Let $A$ be the arithmetic
monodromy group of this cover (i.e. $A$ is the Galois group
of the Galois closure of $F(U)/F(V)$) and $G$ the geometric
monodromy group of the cover (the subgroup of $A$ which
acts trivially on the algebraic closure of $F$). 
Let $H$ be the subgroup of $A$ trivial on $F(U)$ (so $|A:H|=n$). Note
that $A/G$ is cyclic.  Let $xG$ be a generator for $A/G$.
It follows from the Cebotarev density theorem (cf. \cite{gtz})
that:

\begin{theorem}   
$$
|f(U(q_a))| = 1 - \delta(\langle x^a,G\rangle, G, H) + O(q^{a/2}).
$$
\end{theorem}

The special  case where $\delta(\langle x^a,G\rangle, G, H)=0$
gives rise to exceptional covers.  In this case, it is not
difficult to show that $f$ is in fact bijective on rational
points. Of course, this cannot  be the case if $x^a \in G$.
See \cite{FGS}, \cite{GMS} and \cite{gurram} 
for more about exceptional covers.
By \cite{gurstevenson}, any group theoretic solution does give rise
to some cover of curves with the appropriate property.

\section{Derangements in a Coset}\label{coset}

In this section, we present a proof of the Guralnick-Wan result
-- Theorem \ref{gwthm}  that is a bit different than the one given
in \cite{gurwan}.  It is more in the spirit of the proof in \cite{cc} and
an unpublished proof of Marty Isaacs (both for the case $A=G$).

Let $G$ be a normal subgroup of $A$ with $A/G$ generated
by $aG$.  Suppose that $A$ and $G$ both act on the finite set
$X$.
Let $f(g)$ be the number of fixed points of $g$ on $X$.
We note the following well known easy result (cf. \cite{gurwan}).

\begin{lemma}
$\sum_{g \in G} f(ag) =|G|c$, where $c$ is the number of
common $A,G$-orbits on $X$.
\end{lemma}

Now suppose that $A$ and $G$ are both transitive on $X$
(and so in particular $c=1$ in the previous result).
Let $H$ be the stabilizer of a point and set $K=H \cap G$.

Let $\Delta$ denote the derangements in the coset $xG$.
There must be some element in the coset with a fixed point and so we
may assume that $a \in H$.

We split $xG$ into three disjoint sets, $xK$, $\Delta$ and $\Gamma$
(the complement of the union of $xK$ and $\Delta$).

Breaking up the sum into the sum into two pieces, one over $xK$ and
the other the remaining terms, we see that

$$
|G| = \sum_{g \in K} f(ag) + 
\sum_{xg \in \Gamma}  f(ag) \ge c|K| + |G| -|K|- |\Delta|,
$$
where $d$ is the number of common $H,K$ orbits on $X$ (of course, $d \ge 1$).

This yields
$|\Delta| \ge (d-1)|K|$ or $\delta(A,G,X) \ge (d-1)/n$.

If $d=1$, then it is easy to see that $\Delta$ is empty 
(using the fact that the
average number of fixed points is $1$).  So we obtain:

\begin{theorem}  If $(A,G,X)$ is not exceptional, 
then $\delta(A,G,X) \ge 1/n$.
\end{theorem}

If $d \ge 3$, we see that $\delta(A,G,X) \ge 2/n$.  It would be interesting to
characterize those groups where $d=2$ (this includes the case where
$G$ is $2$-transitive) and classify the actions where $\delta(A,G,X) \le 2/n$
(presumably only Frobenius groups and exceptional actions).

\section{Primitive Groups}\label{primitivesec}

As we have seen in the previous section, we cannot hope
to extend Shalev's conjecture to the almost simple case.
There are many more examples in case of affine primitive
groups and also diagonal actions (again related to
exceptionality -- see \cite{GMS} for examples).

In this section we show how one can obtain a weaker result
for primitive groups with no normal abelian subgroup (so
in particular as long as the degree is not a prime power).
The example in the previous section shows that one can
do no better than this theorem.  We do hope to classify
which primitive representations have few derangements.

\begin{theorem}  Let $G$ be a primitive group of degree $n$
and assume that $G$ has no normal abelian subgroup.  Then
there exists a positive constant $\delta$ such that
$\delta(G,X) > \delta/\log  n$.
\end{theorem}

We prove this by reducing to the almost simple case and
then to the simple case.

We first need some auxiliary lemmas.

We will use the following result (which depends on
the classification of finite simple groups --
see  \cite{GMS} for a proof).

\begin{lemma} \label{fixedpoint} If $h$ is an automorphism of a 
finite nonsolvable group $J$, then $C_J(h) \ne 1$.
\end{lemma}

\begin{lemma}  Let $G$ be a transitive permutation group
with a regular nonsolvable  normal subgroup $N$ acting
on $X$.  Then $\delta(G,X) \ge 1/2$.
\end{lemma}

\begin{proof} We can identify $N$ with $X$.
A point stabilizer $H$ is a complement to $N$ and the
action on $X$ is equivalent to the conjugation action
of $H$.  If $h \in H$, then the number of fixed
points is just $|C_N(h)| > 1$ (by the previous result).
Thus every element either has zero
or at least $2$ fixed points.  Since the average number
of fixed points is $1$, this implies that $\delta(G,X) \ge 1/2$.
\end{proof}

We say that $G$ preserves a product structure on $X$ if
$X$ can be identified with $Y \times \ldots \times Y$ ($t > 1$
copies) and $G$ embeds in $S_Y \wr S_t$ in its natural action
($S_Y$ is the symmetric group on $Y$ and each of the $t$ copies
acts on one copy of $Y$, the $S_t$ permutes the coordinates).
In particular, there is a homomorphism $\pi$ from $G$ into $S_t$.
We assume that this image is transitive (which is always the
case if $G$ is primitive on $X$).  Let $G_1$ denote the
preimage of the stabilizer of $1$ in $\pi(G)$.  So $G_1$
acts on $Y$.  If $G$ is primitive, it follows that $G_1$
is as well \cite{aschscott}.

\begin{lemma} $\delta(G,X) \ge \delta(G_1,Y)/t$.
\end{lemma}

\begin{proof} The proportion of elements in $G_1$
is $1/t$.  If $g \in G_1$, then $g$ a derangement
on $Y$ implies that $g$ is a derangement on $X$.
\end{proof}

Note in particular that $\log |X| = t \log |Y|$.

By examining the possibilities of primitive permutation
groups (see \cite{aschscott}) and using the two
previous lemmas, there are only two cases remaining --
$G$ is almost simple or $X$ is of full diagonal type
(we explain this more fully below). Let
$H$ be a point stabilizer.  In particular,
$G$ has a unique minimal normal subgroup $N$ a direct
product of $t$ copies of a nonabelian simple $L$ and either
$t=1$ or we may view $H \cap N \cong L$ as the diagonal
subgroup of $N$ (note that all diagonal subgroups are 
conjugate in $\aut(N)$ so there is no loss of generality
in assuming that $H \cap N$ is the canonical diagonal subgroup
-- alternatively, the arguments below are valid with
$H \cap N$ any diagonal subgroup).

We next handle the diagonal case.

\begin{lemma} Let $G$ be a finite group
with a unique minimal normal subgroup $N=L^{t}$ with $L$ a nonabelian
finite simple group and $t > 1$.  Let $D$ be a diagonal subgroup
of $N$ and assume that $G=NH$ with $H=N_G(D)$.  Then
$\delta(G,G/H) > 1/ \log_2 |G/H|  $.
\end{lemma}

\begin{proof}  Suppose that $g \in G$ has a unique fixed point
on $G/H$.  We claim that $g$ is transitive on the $t$ conjugates
of $L$.  Conjugating by an element of the transitive
subgroup $N$ allows us to assume that $g \in H$.
Since $g$ has a unique fixed point, it is invariant
under $C_G(g)$ and  so $C_G(g) \le H$.
In particular, $C_N(g) \le D$.  This implies the claim
-- for if $g$ leaves invariant
some proper factor $N_1$ of $N$, then $C_{N_1}(g) \ne 1$ (by
Lemma \ref{fixedpoint}) but $N_1 \cap D = 1$.

Now the proportion of elements in $G$ that induce at $t$-cycle
on the $t$ conjugates of $L$ is at most $1 - 1/t$
($1/t$ of the elements normalize $L$).  Thus, the proportion
of elements with a unique fixed point is at most $1-1/t$, whence
at least half the remaining elements must be derangements.
Thus, the proportion of derangements is at least $1/2t$.
Since $|G:H| = |L|^{t-1} \ge 60^{t-1}$, we have $2t < (t-1) \log_2 60
\le \log_2 |G/H|$.

\end{proof}

One can show that in most cases above, one can obtain an estimate
not involving the $\log$ term.  However, if the action of $G$
on the $t$-conjugates of $L$ is cyclic of order $t$ and $t$ does not
divide the order of $L$, then in fact one can do no better than
the previous result (this is another example of exceptionality).

If $G$ is almost simple, then we can apply Theorem \ref{shalevconj}.
Note that this implies the same result for almost simple
groups as long as the socle of $G$ has bounded index (with
perhaps a different constant).  Indeed, in the sequel we actually
prove the result for all Chevalley groups contained in the group
of inner diagonal automorphisms. Since the group of graph
automorphisms always has order at most $6$,  we only need worry about
field automorphisms.  A simple inspection shows that
 the group of the field automorphisms has order
at most $\log_2 n$ where $n$ is the degree of the permutation
representation.  Thus, we have proved our result follows
from Theorem \ref{shalevconj}.  We have proved Theorem \ref{shalevconj}
in the bounded rank case.  As we noted
in the introduction, the complete proof of Theorem \ref{shalevconj}
is contained in \cite{FG1}, \cite{FG2} and \cite{FG3}.

\section{Numbers of conjugacy classes in finite classical groups}
\label{conjugacyclasses}

To conclude this paper we record some upper bounds on the
number of conjugacy classes in the finite classical groups. These are
treated fully in \cite{FG2} where the results are used as a key
ingredient in proving Theorem \ref{shalevconj} and more. 
We mention that
upper bounds on numbers of conjugacy classes are also of interest in
random walks \cite{gluck}, \cite{LiSh}  and
for computation of Fourier transforms on finite groups \cite{MR}.

The bounds we present are of the form $cq^r$ where $r$ is the
rank and $c$ is a small explicit constant. The paper \cite{liepyber} had
previously established the bound $(6q)^r$ and the paper \cite{gluck}
had established bounds such as $cq^{3r}$. Thus our bounds in Theorem
\ref{simpleform} are sharper for classical groups. 
For exceptional groups, one can compute precisely the number of classes
as a monic polynomial in $q$ (see \cite{ssclasses} for some discussion of
this and references).  In even characteristic $O(2n+1,q)$ is
isomorphic to $Sp(2n,q)$ so we omit this case. Here we only state the
results for a specific form of each group.  This gives
bounds for the simple groups using the fact that the number of conjugacy
classes decreases when one takes homomorphic images
and also using  the lemma below to pass to a subgroup or overgroup
of bounded index.  In \cite{FG2}, we actually prove results for more forms
of the groups.

\begin{theorem} \label{simpleform} Let $k(G)$ denote the
number of conjugacy classes of a finite group $G$.
\begin{enumerate}
\item $k(SL(n,q)) \leq \frac{q^n}{q-1} + q^{1+\frac{n}{2}}$.
\item $k(SU(n,q)) \leq 11.5  \left( \frac{q^n}{q+1}+\frac{q+1}{q-1}
q^{n/2+1} \right)$.
\item $k(Sp(2n,q)) \leq 12 q^n$ if $q$ is odd.
\item $k(Sp(2n,q)) \leq 21.4 q^n$ if $q$ is even.
\item $k(O^{\pm}(2n,q)) \leq 29 q^n$ if $q$ is odd.
\item $k(O^{\pm}(2n,q)) \leq 19.5 q^n$ if $q$ is even.
\item $k(SO(2n+1,q)) \leq 7.38 q^n$ if $q$ is odd.
\end{enumerate}
\end{theorem}

Our proof of Theorem \ref{simpleform} uses generating
functions for numbers of conjugacy classes in finite classical groups
\cite{Lu}, \cite{Mac}, \cite{MR}, \cite{W} and is largely inspired
by the proof in  \cite{MR} that $GL(n,q)$ has at most
$q^n$ classes and that $GU(n,q)$ has at most 
$8.26q^n$ conjugacy classes.  However some new ingredients
(combinatorial identities) are required.

Let $k_p(G)$ denote the number of conjugacy classes of
$p'$-elements of $G$.  This is also the number of absolutely
irreducible representations of $G$ in characteristic $p$.  If
$p$ does not divide $G$, then $k_p(G)=k(G)$ the number of conjugacy
classes of $G$ (and also the number of irreducible complex representations).
We also employ the following useful lemma,
which allows us to move between various forms of the finite classical
groups--at least when $|G/H|$ is bounded.  This is proved in
\cite{Ga} for $k(G)$. The modification for $p'$-classes
is straightforward and we omit the proof.

\begin{lemma} \label{Boblemma}  Let $H$ be a subgroup of $G$ with
$G/H$ of order $d$.  Fix a prime $p$.  
Then $k_p(G) \le dk_p(H)$ and $k_p(H) \le d k_p(G) $.
If $H$ is normal in $G$, then $k_p(G) \le k_p(H)k_p(G/H)$.
\end{lemma}

In fact using generating functions it is possible to
understand the asymptotic behavior of the constant $c$ in the bound
$cq^n$ of Theorem \ref{simpleform}.
More precisely, we establish in
\cite{FG2} the following result.

\begin{theorem}
\begin{enumerate}
\item $lim_{n \rightarrow \infty} \frac{k(GL(n,q))}{q^n} =
1$
\item $lim_{n \rightarrow \infty}
\frac{k(GU(n,q))}{q^n} = \prod_{i \geq 1} \frac{1+1/q^i}{1-1/q^i}$
\item $lim_{n \rightarrow \infty} \frac{k(Sp(2n,q))}{q^n} =
\prod_{i=1}^{\infty} \frac{(1+\frac{1}{q^i})^4}{(1-\frac{1}{q^i})}$
if $q$ is odd.
\item $lim_{n \rightarrow \infty} \frac{k(Sp(2n,q))}{q^n} =
\prod_{i=1}^{\infty} \frac{1-1/q^{4i}}{(1-1/q^{4i-2}) (1-1/q^i)^2}$
if $q$ is even.
\item $lim_{n \rightarrow \infty} \frac{k(O^{\pm}(2n,q))}{q^{n}} =
\frac{1}{4 \prod_{i=1}^{\infty} (1-1/q^{i})} (\prod_{i=1}^{\infty}
(1+1/q^{i-1/2})^4 +\prod_{i=1}^{\infty} (1-1/q^{i-1/2})^4)$ if $q$ is odd.
\item $lim_{n \rightarrow \infty} \frac{k(O^{\pm}(2n,q))}{q^n} = \frac{1}{2}
\frac{\prod_{i=0}^{\infty} (1-1/q^{2i+2})
(1+1/q^{2i+1})^2}{\prod_{i=1}^{\infty} (1-1/q^i)^2}$ if $q$ is even.
\item $lim_{n \rightarrow \infty} \frac{k(SO(2n+1,q))}{q^n} =
\prod_{i=1}^{\infty} \frac{(1-1/q^{4i})^2}{(1-1/q^i)^3
(1-1/q^{4i-2})^2}$ if $q$ is odd.
\end{enumerate}
\end{theorem}

\end{document}